\documentclass[10pt,oneside,english]{amsart}
\usepackage[T1]{fontenc}
\usepackage[latin9]{inputenc}
\usepackage{geometry}
\usepackage{babel}
\usepackage{amstext}
\usepackage{amsthm}
\usepackage{amssymb}
\usepackage[unicode=true,pdfusetitle,
 bookmarks=true,bookmarksnumbered=false,bookmarksopen=false,
 breaklinks=false,pdfborder={0 0 1},backref=false,colorlinks=false]
 {hyperref}

\makeatletter
\numberwithin{equation}{section}
\numberwithin{figure}{section}
\theoremstyle{plain}
\newtheorem{thm}{\protect\theoremname}
  \theoremstyle{remark}
  \newtheorem{rem}[thm]{\protect\remarkname}
  \theoremstyle{remark}
  \newtheorem*{acknowledgement*}{\protect\acknowledgementname}

\@ifundefined{date}{}{\date{}}

\makeatother

  \providecommand{\acknowledgementname}{Acknowledgement}
  \providecommand{\remarkname}{Remark}
\providecommand{\theoremname}{Theorem}

\begin{document}

\title{Regularity bounds for a Gevrey criterion in a kernel-based regularization
of the Cauchy problem of elliptic equations}

\author{Vo Anh Khoa}

\address{Mathematics and Computer Science Division, Gran Sasso Science Institute,
L'Aquila, Italy.}

\email{khoa.vo@gssi.infn.it, vakhoa.hcmus@gmail.com}

\author{Tran The Hung}

\address{Faculty of Applied Physics and Mathematics, Gdasnk University of
Technology, Gdasnk, Poland.}

\email{thehung.tran@mathmods.eu}

\keywords{Gevrey criterion, Regularity bounds, Kernel-based regularization,
Cauchy problem, Elliptic equations, Inequalities}

\subjclass[2000]{26D15, 35J91, 47A52, 46E35}
\begin{abstract}
This Note derives regularity bounds for a Gevrey criterion when the
Cauchy problem of elliptic equations is solved by regularization.
When utilizing the regularization, one knows that checking such criterion
is basically problematic, albeit its importance to engineering circumstances.
Therefore, coping with that impediment helps us improve the use of
some regularization methods in real-world applications. This work
also consider the presence of the power-law nonlinearities.
\end{abstract}

\maketitle

\section{Background}

Let us consider the Cauchy problem of semi-linear elliptic equations,
as follows:
\[
\frac{d^{2}\mathbf{u}\left(x,y\right)}{dx^{2}}=\mathcal{A}_{y}\mathbf{u}\left(x,y\right)+f\left(x,y,\mathbf{u}\left(x,y\right)\right),\quad\left(x,y\right)\in\Omega:=\Omega_{x}\times\Omega_{y},
\]
associated with the zero Dirichlet boundary condition in $y$ and
the initial data and nonhomogeneous initial velocity given by
\[
\mathbf{u}\left(0,y\right)=\mathbf{u}_{0}\left(y\right),\quad\frac{d\mathbf{u}\left(0,y\right)}{dx}=\mathbf{u}_{1}\left(y\right),\quad y\in\Omega_{y}.
\]

Here, $\mathbf{u}:\Omega_{x}\to L^{2}\left(\Omega_{y}\right)$ denotes
the distribution of a body where $\Omega_{x}:=\left(0,a\right)\subset\mathbb{R}$
and $\Omega_{y}\subset\mathbb{R}^{n}$ are open, bounded and connected
domains for $n\ge2$ and $a>0$ with a smooth boundary, and $\mathcal{A}_{y}$
is the linear second-order differential operator with variable coefficients
depending on $y$ only:
\[
\mathcal{A}_{y}\mathbf{u}\left(x,y\right)=\sum_{i,j=1}^{n}\frac{\partial}{\partial y_{i}}\left(d_{i,j}\left(y\right)\frac{\partial\mathbf{u}\left(x,y\right)}{\partial y_{j}}\right)+d\left(y\right)u\left(y\right).
\]

The basic requirement for the coefficients $d_{i,j}\left(y\right)$
and $d\left(y\right)$ is that $\mathcal{A}_{y}$ is a positive, self-adjoint
operator in the Hilbert space $L^{2}\left(\Omega_{y}\right)$. Consequently,
there exists an orthonormal basis of $L^{2}\left(\Omega_{y}\right)$,
denoted by $\left\{ \phi_{p}\right\} _{p\in\mathbb{N}^{*}}$, satisfying
\begin{equation}
\phi_{p}\in H_{0}^{1}\left(\Omega_{y}\right)\cap C^{\infty}\left(\overline{\Omega_{y}}\right),\quad\mathcal{A}_{y}\phi_{p}\left(y\right)=\lambda_{p}\phi_{p}\left(y\right)\;\text{for}\;y\in\Omega_{y},\label{eq:eigen1}
\end{equation}
and the corresponding discrete spectrum $\left\{ \lambda_{p}\right\} _{p\in\mathbb{N}^{*}}$
satisfies
\begin{equation}
0<\lambda_{1}\le\lambda_{2}\le...\lim_{p\to\infty}\lambda_{p}=\infty.\label{eigen2}
\end{equation}

As a direct example, we take $\mathcal{A}_{y}=-\Delta_{y}$ with an
open parallelepiped $\Omega_{y}=\left(0,a_{1}\right)\times...\times\left(0,a_{n}\right)\subset\mathbb{R}^{n}$.
For each $p\in\mathbb{N}^{*}$, it thus gives us that 
\begin{equation}
\phi_{p}\left(y_{1},...,y_{n}\right)=\prod_{j=1}^{n}\sqrt{\frac{2}{a_{j}}}\sin\left(\frac{\pi n_{j}}{a_{j}}y_{j}\right),\quad\lambda_{p}=\sum_{j=1}^{n}\left(\frac{\pi n_{j}}{a_{j}}\right)^{2},\quad n_{j}\in\mathbb{N},j\in\left\{ 1,...,n\right\} ,\label{eq:eigennn}
\end{equation}
which fulfill (\ref{eq:eigen1}) and (\ref{eigen2}), respectively.

It is worth mentioning that this kind of problems is widely known
and esteemed. Essentially, it includes the elliptic sine\textendash Gordon
equations in superconductivity, the Lane\textendash Emden\textendash Fowler
type system arising in molecular biology and the Helmholtz equation
together with its modified versions. For ease of presentation, we
refer the above-mentioned Cauchy problem as Problem $\left(P\right)$.
In this Note, we are interested in the mild solution for $\left(P\right)$
where solving it can be found in \cite{TTK15}, e.g. and then with
$\left(\mathbf{u}_{0},\mathbf{u}_{1}\right)\in L^{2}\left(\Omega_{y}\right)\times L^{2}\left(\Omega_{y}\right)$
we obtain
\begin{align}
\mathbf{u}\left(x,y\right) & =\sum_{p=1}^{\infty}\left[\cosh\left(\sqrt{\lambda_{p}}x\right)\left\langle \mathbf{u}_{0},\phi_{p}\right\rangle +\frac{\sinh\left(\sqrt{\lambda_{p}}x\right)}{\sqrt{\lambda_{p}}}\left\langle \mathbf{u}_{1},\phi_{p}\right\rangle \right.\nonumber \\
 & \left.+\int_{0}^{x}\frac{\sinh\left(\sqrt{\lambda_{p}}\left(x-\xi\right)\right)}{\sqrt{\lambda_{p}}}\left\langle f\left(\xi,\cdot,\mathbf{u}\left(\xi,\cdot\right)\right),\phi_{p}\right\rangle d\xi\right]\phi_{p}\left(y\right),\quad\left(x,y\right)\in\Omega,\label{eq:1.3}
\end{align}
where $\left\langle \cdot,\cdot\right\rangle $ denotes the inner
product in $L^{2}\left(\Omega_{y}\right)$.

Hereby, it is not difficult to see from (\ref{eq:1.3}) that $\cosh\left(\sqrt{\lambda_{p}}x\right)$
and $\sinh\left(\sqrt{\lambda_{p}}x\right)/\sqrt{\lambda_{p}}$ are
all unbounded terms. As a result, their catastrophic growth (as $p\to\infty$)
ruins any computations on the solution $\mathbf{u}\left(x,y\right)$.
In addition, one usually meets the measurement in practice, i.e. we
need to assume the presence of an approximation $\left(\mathbf{u}_{0}^{\varepsilon},\mathbf{u}_{1}^{\varepsilon}\right)\in L^{2}\left(\Omega_{y}\right)\times L^{2}\left(\Omega_{y}\right)$
that satifies
\begin{equation}
\left\Vert \mathbf{u}_{0}^{\varepsilon}-\mathbf{u}_{0}\right\Vert _{L^{2}\left(\Omega_{y}\right)}+\left\Vert \mathbf{u}_{1}^{\varepsilon}-\mathbf{u}_{1}\right\Vert _{L^{2}\left(\Omega_{y}\right)}\le\varepsilon,\label{eq:epsilon}
\end{equation}
in which the constant $\varepsilon>0$ represents the upper bound
of the noise level in $L^{2}\left(\Omega_{y}\right)$.

In order to overcome the Hadamard-instability for this type of problems,
some regularization methods have been proposed: the quasi-reversibility
method \cite{QFL08}, the quasi-boundary value method \cite{HDL09}
and the truncation method \cite{TTQ10}. We notice herein that when
using the kernel-based regularization, the Gevrey criterion is faced.
In particular, we consider the following regularized solution
\begin{align}
\mathbf{u}^{\varepsilon}\left(x,y\right) & =\sum_{p=1}^{\infty}\left[\cosh^{\varepsilon}\left(\sqrt{\lambda_{p}}x\right)\left\langle \mathbf{u}_{0},\phi_{p}\right\rangle +\frac{\sinh^{\varepsilon}\left(\sqrt{\lambda_{p}}x\right)}{\sqrt{\lambda_{p}}}\left\langle \mathbf{u}_{1},\phi_{p}\right\rangle \right.\nonumber \\
 & \left.+\int_{0}^{x}\frac{\sinh^{\varepsilon}\left(\sqrt{\lambda_{p}}\left(x-\xi\right)\right)}{\sqrt{\lambda_{p}}}\left\langle f\left(\xi,\cdot,\mathbf{u}^{\varepsilon}\left(\xi,\cdot\right)\right),\phi_{p}\right\rangle d\xi\right]\phi_{p}\left(y\right),\quad\left(x,y\right)\in\Omega,\label{eq:regu}
\end{align}
where for each $\varepsilon>0$ the terms $\cosh^{\varepsilon}\left(\sqrt{\lambda_{p}}x\right)$
and $\sinh^{\varepsilon}\left(\sqrt{\lambda_{p}}x\right)/\sqrt{\lambda_{p}}$
can be bounded from above. This also leads to the conditional stability
estimate for the regularized solution. In recent works, they are of
the form
\[
\cosh^{\varepsilon}\left(\sqrt{\lambda_{p}}x\right):=\Psi_{p,k}^{\beta}\left(x\right)+\frac{e^{-\sqrt{\lambda_{p}}x}}{2},\quad\sinh^{\varepsilon}\left(\sqrt{\lambda_{p}}x\right)=\Psi_{p,k}^{\beta}\left(x\right)-\frac{e^{-\sqrt{\lambda_{p}}x}}{2},
\]
in which the general kernel $\Psi_{p,k}^{\beta}\left(x\right):\overline{\Omega_{x}}\to\mathbb{R}_{+}$
is provided by
\begin{equation}
\Psi_{p,k}^{\beta}\left(x\right):=\frac{e^{-\sqrt{\lambda_{p}}\left(a-x\right)}}{2\beta\sqrt{\lambda_{p}^{k}}+2e^{-\sqrt{\lambda_{p}}a}},\quad p\in\mathbb{N}^{*},\beta:=\beta\left(\varepsilon\right)\in\left(0,1\right),\label{eq:kernel}
\end{equation}
with $k\ge1$ inspired from \cite{KTDT17} and \cite{TTTK15}, and
$k=0$ postulated in \cite{TTK15}.

When doing so, the Gevrey criterion for convergence is known as the
\emph{a priori} information on the exact solution under the Gevrey\footnote{Here, we employ this terminology from Cao et al. \cite{CRT99}.}
classes defined by
\[
\mathbb{G}_{\nu}^{s}:=\left\{ \mathbf{v}\in L^{2}\left(\Omega_{y}\right):\sum_{p=1}^{\infty}\lambda_{p}^{\nu}e^{2s\sqrt{\lambda_{p}}}\left|\left\langle \mathbf{v},\phi_{p}\right\rangle \right|^{2}<\infty\right\} ,\quad\nu\ge0,s>0,
\]
endowed with the norm
\[
\left\Vert \mathbf{v}\right\Vert _{\mathbb{G}_{\nu}^{s}}^{2}=\sum_{p=1}^{\infty}\lambda_{p}^{\nu}e^{2s\sqrt{\lambda_{p}}}\left|\left\langle \mathbf{v},\phi_{p}\right\rangle \right|^{2}<\infty.
\]

Return to our concern, from \cite[Theorem 7]{KTDT17} and \cite[Theorem 2]{TTTK15}
it requires that
\begin{equation}
\mathbf{u}\in C\left(\overline{\Omega_{x}};\mathbb{G}_{\nu_{1}}^{s_{1}}\right)\;\text{and}\;\frac{d\mathbf{u}}{dx}\in C\left(\overline{\Omega_{x}};\mathbb{G}_{\nu_{2}}^{s_{2}}\right),\label{eq:1.4}
\end{equation}
with $\nu_{1}\ge a$, $\nu_{2}\ge a$, $s_{1}=k$ and $s_{2}=k-1$,
whilst in \cite{TTK15} we assume that
\begin{equation}
\mathbf{u}\in C\left(\overline{\Omega_{x}};\mathbb{G}_{0}^{s_{1}}\right)\;\text{and}\;\frac{d\mathbf{u}}{dx}\in C\left(\overline{\Omega_{x}};\mathbb{G}_{0}^{s_{2}}\right),\label{eq:1.5}
\end{equation}
with $s_{1}$ and $s_{2}$ being the same as above.

At present, we observe that the assumptions (\ref{eq:1.4}) and (\ref{eq:1.5})
are very hard to check if one wants to utilize this type of regularization
methods and from those works mentioned above, they merely consider
this information when $f\equiv0$. Due to those reasons, this Note
is to explore a natural upper bound for such criterion for $k\in\mathbb{N}^{*}$
and due to the similarity, we focus on the assumption (\ref{eq:1.4})
in the next section. Our main result thus lies in Theorem \ref{thm:Consider-the-general}.

\section{Derivation of regularity bounds}

For simplicity, the forcing function $f\left(x,y,\mathbf{u}\left(x,y\right)\right)=\mathbf{f}\left(\mathbf{u}\left(x,y\right)\right)+\mathbf{F}\left(x,y\right)$
is concentrated with $\mathbf{f}\left(0\right)\equiv0$. Furthermore,
we assume a modulus of continuity $\omega:\left[0,\infty\right]\to\left[0,\infty\right]$
on $\mathbf{f}:\mathbb{R}\to\mathbb{R}$, i.e.
\begin{equation}
\left|\mathbf{f}\left(\mathbf{u}\right)-\mathbf{f}\left(\mathbf{v}\right)\right|\le\omega\left(\left|\mathbf{u}-\mathbf{v}\right|\right)\quad\text{for all}\;\mathbf{u},\mathbf{v}\in\mathbb{R}.\label{eq:conti}
\end{equation}

In this part, we mostly take into consideration the modulus $\omega\left(\mathbf{u}\right):=L\mathbf{u}$
which indicates the globally $L$-Lipschitz function, whilst the H\"older-type
continuity $\omega\left(\mathbf{u}\right):=L\mathbf{u}^{\alpha}$,
$\alpha\ge1$ resembling the power-law nonlinearities (e.g. logistic
and von Bertalanffy) shall be investigated in a few words as a consequence.

From here on, we recall from the proofs of \cite[Theorem 7]{KTDT17}
and \cite[Theorem 2]{TTTK15} the actual assumption that leads to
(\ref{eq:1.4}). It has the following form:
\begin{equation}
\mathbf{A}:=\sup_{x\in\overline{\Omega_{x}}}\sum_{p=1}^{\infty}\lambda_{p}^{k}e^{2\sqrt{\lambda_{p}}\left(a-x\right)}\left(\left\langle \mathbf{u}\left(x,\cdot\right),\phi_{p}\right\rangle +\frac{\left\langle \mathbf{u}_{x}\left(x,\cdot\right),\phi_{p}\right\rangle }{\sqrt{\lambda_{p}}}\right)^{2}<\infty.\label{eq:cond}
\end{equation}

From (\ref{eq:1.3}), we take the derivative of $\mathbf{u}\left(x,y\right)$
with respect to $x$ and obtain that
\begin{align}
\frac{\left\langle \mathbf{u}_{x}\left(x,\cdot\right),\phi_{p}\right\rangle }{\sqrt{\lambda_{p}}} & =\sinh\left(\sqrt{\lambda_{p}}x\right)\left\langle \mathbf{u}_{0},\phi_{p}\right\rangle +\frac{\cosh\left(\sqrt{\lambda_{p}}x\right)}{\sqrt{\lambda_{p}}}\left\langle \mathbf{u}_{1},\phi_{p}\right\rangle \nonumber \\
 & +\int_{0}^{x}\frac{\cosh\left(\sqrt{\lambda_{p}}\left(x-\xi\right)\right)}{\sqrt{\lambda_{p}}}\left\langle f\left(\xi,\cdot,\mathbf{u}\left(\xi,\cdot\right)\right),\phi_{p}\right\rangle d\xi,\quad p\in\mathbb{N}^{*}.\label{eq:2.1}
\end{align}

Therefore, combining (\ref{eq:2.1}) with (\ref{eq:1.3}), we arrive
at
\begin{align}
\left\langle \mathbf{u}\left(x,\cdot\right),\phi_{p}\right\rangle +\frac{\left\langle \mathbf{u}_{x}\left(x,\cdot\right),\phi_{p}\right\rangle }{\sqrt{\lambda_{p}}} & =e^{\sqrt{\lambda_{p}}x}\left(\left\langle \mathbf{u}_{0},\phi_{p}\right\rangle +\frac{1}{\lambda_{p}}\left\langle \mathbf{u}_{1},\phi_{p}\right\rangle \right.\nonumber \\
 & \left.+\int_{0}^{x}\frac{e^{-\sqrt{\lambda_{p}}\xi}}{\sqrt{\lambda_{p}}}\left\langle f\left(\xi,\cdot,\mathbf{u}\left(\xi,\cdot\right)\right),\phi_{p}\right\rangle d\xi\right),\quad p\in\mathbb{N}^{*}.\label{eq:2.2}
\end{align}

Taking now $x=a$ in (\ref{eq:2.2}), we can write that
\begin{align}
\left\langle \mathbf{u}\left(a,\cdot\right),\phi_{p}\right\rangle +\frac{\left\langle \mathbf{u}_{x}\left(a,\cdot\right),\phi_{p}\right\rangle }{\sqrt{\lambda_{p}}} & =e^{\sqrt{\lambda_{p}}a}\left(\left\langle \mathbf{u}_{0},\phi_{p}\right\rangle +\frac{1}{\lambda_{p}}\left\langle \mathbf{u}_{1},\phi_{p}\right\rangle \right.\nonumber \\
 & \left.+\int_{0}^{a}\frac{e^{-\sqrt{\lambda_{p}}\xi}}{\sqrt{\lambda_{p}}}\left\langle f\left(\xi,\cdot,\mathbf{u}\left(\xi,\cdot\right)\right),\phi_{p}\right\rangle d\xi\right),\quad p\in\mathbb{N}^{*}.\label{eq:2.2-1}
\end{align}

Henceforward, combining (\ref{eq:2.2}) and (\ref{eq:2.2-1}) we gain
the following equality after some arrangements
\begin{align*}
e^{\sqrt{\lambda_{p}}\left(a-x\right)}\left(\left\langle \mathbf{u}\left(x,\cdot\right),\phi_{p}\right\rangle +\frac{\left\langle \mathbf{u}_{x}\left(x,\cdot\right),\phi_{p}\right\rangle }{\sqrt{\lambda_{p}}}\right) & =\left\langle \mathbf{u}\left(a,\cdot\right),\phi_{p}\right\rangle +\frac{\left\langle \mathbf{u}_{x}\left(a,\cdot\right),\phi_{p}\right\rangle }{\sqrt{\lambda_{p}}}\\
 & -\int_{x}^{a}\frac{e^{\sqrt{\lambda_{p}}\left(a-\xi\right)}}{\sqrt{\lambda_{p}}}\left\langle f\left(\xi,\cdot,\mathbf{u}\left(\xi,\cdot\right)\right),\phi_{p}\right\rangle d\xi,\quad p\in\mathbb{N}^{*}.
\end{align*}

Hereby, we bound $\mathbf{A}$ from above by
\begin{align}
\mathbf{A} & \le\sup_{x\in\overline{\Omega_{x}}}\sum_{p=1}^{\infty}\lambda_{p}^{k}\left(\left\langle \mathbf{u}\left(a,\cdot\right),\phi_{p}\right\rangle +\frac{\left\langle \mathbf{u}_{x}\left(a,\cdot\right),\phi_{p}\right\rangle }{\sqrt{\lambda_{p}}}-\int_{x}^{a}\frac{e^{\sqrt{\lambda_{p}}\left(a-\xi\right)}}{\sqrt{\lambda_{p}}}\left\langle f\left(\xi,\cdot,\mathbf{u}\left(\xi,\cdot\right)\right),\phi_{p}\right\rangle d\xi\right)^{2}\nonumber \\
 & \le3\sup_{x\in\overline{\Omega_{x}}}\sum_{p=1}^{\infty}\left[\lambda_{p}^{k}\left|\left\langle \mathbf{u}\left(a,\cdot\right),\phi_{p}\right\rangle \right|^{2}+\lambda_{p}^{k-1}\left(\left|\left\langle \mathbf{u}_{x}\left(a,\cdot\right),\phi_{p}\right\rangle \right|^{2}+\left|\int_{x}^{a}e^{\sqrt{\lambda_{p}}\left(a-\xi\right)}\left\langle f\left(\xi,\cdot,\mathbf{u}\left(\xi,\cdot\right)\right),\phi_{p}\right\rangle d\xi\right|^{2}\right)\right],\label{2.6}
\end{align}
where we use the elementary inequality $\left(a_{1}+a_{2}+a_{3}\right)^{2}\le3\left(a_{1}^{2}+a_{2}^{2}+a_{3}^{2}\right)$.

At this point, we observe the norm of the Hilbert space $H^{r}\left(\Omega_{y}\right)$
with $r\in\mathbb{N}$, which can be naturally defined in terms of
Fourier series whose coefficients that decay rapidly; namely
\[
H^{r}\left(\Omega_{y}\right):=\left\{ \mathbf{v}\in L^{2}\left(\Omega_{y}\right):\left\Vert \mathbf{v}\right\Vert _{H^{r}\left(\Omega_{y}\right)}<\infty\right\} ,
\]
equipped with the norm
\[
\left\Vert \mathbf{v}\right\Vert _{H^{r}\left(\Omega_{y}\right)}^{2}=\sum_{p=1}^{\infty}\left(1+\lambda_{p}\right)^{r}\left|\left\langle \mathbf{v},\phi_{p}\right\rangle \right|^{2}.
\]

It then sufficient to bound the criterion $\mathbf{A}$ from above.
Indeed, we get that
\begin{align}
\frac{1}{3}\mathbf{A} & \le\left\Vert \mathbf{u}\left(a,\cdot\right)\right\Vert _{H^{k}\left(\Omega_{y}\right)}^{2}+\left\Vert \mathbf{u}_{x}\left(a,\cdot\right)\right\Vert _{H^{k-1}\left(\Omega_{y}\right)}^{2}+\sup_{x\in\overline{\Omega_{x}}}\sum_{p=1}^{\infty}\lambda_{p}^{k-1}\left|\int_{x}^{a}e^{\sqrt{\lambda_{p}}\left(a-\xi\right)}\left\langle f\left(\xi,\cdot,\mathbf{u}\left(\xi,\cdot\right)\right),\phi_{p}\right\rangle d\xi\right|^{2}\nonumber \\
 & \le\left\Vert \mathbf{u}\left(a,\cdot\right)\right\Vert _{H^{k}\left(\Omega_{y}\right)}^{2}+\left\Vert \mathbf{u}_{x}\left(a,\cdot\right)\right\Vert _{H^{k-1}\left(\Omega_{y}\right)}^{2}+\sup_{x\in\overline{\Omega_{x}}}\sum_{p=1}^{\infty}\lambda_{p}^{k-1}\int_{x}^{a}e^{2\sqrt{\lambda_{p}}\left(a-\xi\right)}d\xi\int_{x}^{a}\left|\left\langle f\left(\xi,\cdot,\mathbf{u}\left(\xi,\cdot\right)\right),\phi_{p}\right\rangle \right|^{2}d\xi\nonumber \\
 & \le\left\Vert \mathbf{u}\left(a,\cdot\right)\right\Vert _{H^{k}\left(\Omega_{y}\right)}^{2}+\left\Vert \mathbf{u}_{x}\left(a,\cdot\right)\right\Vert _{H^{k-1}\left(\Omega_{y}\right)}^{2}+\frac{1}{2}\sum_{p=1}^{\infty}\lambda_{p}^{k-\frac{3}{2}}\left(e^{2\sqrt{\lambda_{p}}a}-1\right)\int_{0}^{a}\left|\left\langle f\left(\xi,\cdot,\mathbf{u}\left(\xi,\cdot\right)\right),\phi_{p}\right\rangle \right|^{2}d\xi,\label{2.5}
\end{align}
in which we apply the fundamental inequalities $\lambda_{p}^{k}\le\left(1+\lambda_{p}\right)^{k}$,
$\lambda_{p}^{k-1}\le\left(1+\lambda_{p}\right)^{k-1}$ for all $p\in\mathbb{N}^{*}$,
$k\in\mathbb{N}^{*}$ in combination with the H\"older inequality.

Consequently, (\ref{2.5}) yields
\begin{equation}
\mathbf{A}\le3\left\Vert \mathbf{u}\left(a,\cdot\right)\right\Vert _{H^{k}\left(\Omega_{y}\right)}^{2}+3\left\Vert \mathbf{u}_{x}\left(a,\cdot\right)\right\Vert _{H^{k-1}\left(\Omega_{y}\right)}^{2}+\begin{cases}
\frac{3}{2}\left\Vert f\left(\mathbf{u}\right)\right\Vert _{L^{1}\left(\Omega_{x};\mathbb{G}_{k-\frac{3}{2}}^{a}\right)}^{2}, & k>1,\\
\frac{3}{2}\lambda_{1}^{k-\frac{3}{2}}\left\Vert f\left(\mathbf{u}\right)\right\Vert _{L^{1}\left(\Omega_{x};\mathbb{G}_{0}^{a}\right)}^{2}, & k=1.
\end{cases}\label{eq:2.7}
\end{equation}

Accordingly, the estimate (\ref{eq:2.7}) completes a general regularity
bound for the Gevrey-type criterion $\mathbf{A}$ defined in (\ref{eq:cond}).
In other words, the assumption that we have constructed facilitates
very much the previously used information (\ref{eq:1.4}) since the
Gevrey class just imposes on the forcing function $f$. In the context
of reconstructing the temperature of a body from interior measurement
in linear cases ($\mathbf{f}\equiv0$), we only need to verify the
distribution and its velocity on the surface $x=a$ in $H^{k}\left(\Omega_{y}\right)$
and $H^{k-1}\left(\Omega_{y}\right)$, respectively, together with
the source function $\mathbf{F}$ that substitutes $f$ in (\ref{eq:2.7}).
Furthermore, if $f\equiv0$, one obtains the following equivalence
relation
\[
\left\Vert \mathbf{u}\right\Vert _{C\left(\overline{\Omega_{x}};H^{k}\left(\Omega_{y}\right)\right)}^{2}+\left\Vert \mathbf{u}_{x}\right\Vert _{C\left(\overline{\Omega_{x}};H^{k-1}\left(\Omega_{y}\right)\right)}^{2}\le\left\Vert \mathbf{u}\right\Vert _{C\left(\overline{\Omega_{x}};\mathbb{G}_{k}^{a}\right)}^{2}+\left\Vert \mathbf{u}_{x}\right\Vert _{C\left(\overline{\Omega_{x}};\mathbb{G}_{k-1}^{a}\right)}^{2}\le3\left(\left\Vert \mathbf{u}\right\Vert _{C\left(\overline{\Omega_{x}};H^{k}\left(\Omega_{y}\right)\right)}^{2}+\left\Vert \mathbf{u}_{x}\right\Vert _{C\left(\overline{\Omega_{x}};H^{k-1}\left(\Omega_{y}\right)\right)}^{2}\right).
\]

The regularity bound (\ref{eq:2.7}) is very helpful but we can derive
a more rigorous bound by considering the forcing term $f$. Suppose
$\mathbf{F}\equiv0$, it is straightforward to deduce from (\ref{eq:conti})
that
\[
\mathbf{A}\le3\left(\left\Vert \mathbf{u}\left(a,\cdot\right)\right\Vert _{H^{k}\left(\Omega_{y}\right)}^{2}+\left\Vert \mathbf{u}_{x}\left(a,\cdot\right)\right\Vert _{H^{k-1}\left(\Omega_{y}\right)}^{2}\right)+\begin{cases}
\frac{3L^{2}}{2}\left\Vert \mathbf{u}\right\Vert _{L^{1}\left(\Omega_{x};\mathbb{G}_{k-\frac{3}{2}}^{a}\right)}^{2}, & k>1,\\
\frac{3}{2}\lambda_{1}^{k-\frac{3}{2}}L^{2}\left\Vert \mathbf{u}\right\Vert _{L^{1}\left(\Omega_{x};\mathbb{G}_{0}^{a}\right)}^{2}, & k=1.
\end{cases},
\]
if the modulus $\omega$ is the globally $L$-Lipschitz function (e.g.
$\mathbf{f}\left(\mathbf{u}\right)=\sin\left(\mathbf{u}\right)$ with
$L=1$ and $\mathbf{f}\left(\mathbf{u}\right)=\mathbf{u}\left(1+\mathbf{u}^{2}\right)^{-1}$
with $L=25/16$). This means that we can assume $\mathbf{u}\in C^{1}\left(\overline{\Omega_{x}};H^{k}\left(\Omega_{y}\right)\right)\cap L^{1}\left(\Omega_{x};\mathbb{G}_{k-\frac{3}{2}}^{a}\right)$
if $k\ge2$ and $\mathbf{u}\in C^{1}\left(\overline{\Omega_{x}};H^{k}\left(\Omega_{y}\right)\right)\cap L^{1}\left(\Omega_{x};\mathbb{G}_{0}^{a}\right)$
if $k=1$. On the other side, if $\omega\left(\mathbf{u}\right)=L\mathbf{u}^{\alpha}$
with $\alpha\ge1$ and we know that $\mathbf{u}$ is positive and
bounded, then one can prove $\omega$ is still globally Lipschitz.

All in all, we now state the following theorem.
\begin{thm}
\label{thm:Consider-the-general}Consider the general kernel-based
regularization with kernel $\Psi_{p,k}^{\beta}$ defined in (\ref{eq:kernel})
in accordance with $k\in\mathbb{N}^{*}$ and $\beta:=\beta\left(\varepsilon\right)\in\left(0,1\right)$.
Then, the Gevrey-type criterion (\ref{eq:cond}) on the exact solution
$\mathbf{u}$ of Problem $\left(P\right)$, which are required for
the convergence rate of the regularized solution $\mathbf{u}^{\varepsilon}$
defined in (\ref{eq:regu}) are accepted by the regularity bound
\[
\mathbf{A}\le3\left(\left\Vert \mathbf{u}\left(a,\cdot\right)\right\Vert _{H^{k}\left(\Omega_{y}\right)}^{2}+\left\Vert \mathbf{u}_{x}\left(a,\cdot\right)\right\Vert _{H^{k-1}\left(\Omega_{y}\right)}^{2}\right)+\begin{cases}
\frac{3}{2}\left\Vert f\left(\mathbf{u}\right)\right\Vert _{L^{1}\left(\Omega_{x};\mathbb{G}_{k-\frac{3}{2}}^{a}\right)}^{2}, & k\ge2,\\
\frac{3}{2}\lambda_{1}^{k-\frac{3}{2}}\left\Vert f\left(\mathbf{u}\right)\right\Vert _{L^{1}\left(\Omega_{x};\mathbb{G}_{0}^{a}\right)}^{2}, & k=1.
\end{cases}
\]
Furthermore, if consider $f\left(x,y,\mathbf{u}\left(x,y\right)\right)=\mathbf{f}\left(\mathbf{u}\left(x,y\right)\right)$
with $\mathbf{f}\left(0\right)\equiv0$ and the modulus of continuity
$\omega$ satisfying (\ref{eq:conti}) is globally Lipschitz, the
criterion becomes $\mathbf{u}\in C^{1}\left(\overline{\Omega_{x}};H^{k}\left(\Omega_{y}\right)\right)\cap L^{1}\left(\Omega_{x};\mathbb{G}_{k-\frac{3}{2}}^{a}\right)$
if $k\ge2$ and $\mathbf{u}\in C^{1}\left(\overline{\Omega_{x}};H^{k}\left(\Omega_{y}\right)\right)\cap L^{1}\left(\Omega_{x};\mathbb{G}_{0}^{a}\right)$
if $k=1$.
\end{thm}
\begin{rem}
We define another criterion $\mathbf{A}^{\gamma}$ with an index $\gamma>0$,
provided by
\[
\mathbf{A}^{\gamma}:=\sup_{x\in\overline{\Omega_{x}}}\sum_{p=1}^{\infty}\lambda_{p}^{k}e^{2\gamma\sqrt{\lambda_{p}}\left(a-x\right)}\left(\left\langle \mathbf{u}\left(x,\cdot\right),\phi_{p}\right\rangle +\frac{\left\langle \mathbf{u}_{x}\left(x,\cdot\right),\phi_{p}\right\rangle }{\sqrt{\lambda_{p}}}\right)^{2\gamma}<\infty.
\]

Obviously, this criterion is considered as a special case of $\mathbf{A}$.
Consider the case $\mathbf{f}\equiv0$, we proceed the same way as
estimated in (\ref{2.6}) and (\ref{2.5}). We know that there always
exists a positive constant $C>0$ such that for any countably infinite
set $\left\{ a_{p}\right\} _{p\in\mathbb{N}^{*}}$, whose elements
are all nonnegative, satisfying ${\displaystyle \sum_{p=1}^{\infty}}a_{p}<\infty$,
the following inequality holds
\[
\sum_{p=1}^{\infty}a_{p}^{\gamma}\le C\left(\sum_{p=1}^{\infty}a_{p}\right)^{\gamma}\quad\text{for}\;\gamma>0.
\]

Therefore, it enables us to estimate $\mathbf{A}^{\gamma}$ from above
by
\begin{align}
\frac{1}{2^{\gamma}}\mathbf{A}^{\gamma} & \le\sum_{p=1}^{\infty}\left(\lambda_{p}^{\frac{k}{\gamma}}\left|\left\langle \mathbf{u}\left(a,\cdot\right),\phi_{p}\right\rangle \right|^{2}+\lambda_{p}^{\frac{k}{\gamma}-1}\left|\left\langle \mathbf{u}_{x}\left(a,\cdot\right),\phi_{p}\right\rangle \right|^{2}\right)^{\gamma}\nonumber \\
 & \le C\left(\sum_{p=1}^{\infty}\lambda_{p}^{\frac{k}{\gamma}}\left|\left\langle \mathbf{u}\left(a,\cdot\right),\phi_{p}\right\rangle \right|^{2}+\sum_{p=1}^{\infty}\lambda_{p}^{\frac{k}{\gamma}-1}\left|\left\langle \mathbf{u}_{x}\left(a,\cdot\right),\phi_{p}\right\rangle \right|^{2}\right)^{\gamma}.\label{eq:2.9}
\end{align}

At present, we need to argue the relation between $k$ and $\gamma$.
In fact, if $k>\gamma$ then using the inequality
\[
\left(a+b\right)^{\gamma}\le\max\left\{ 2^{\gamma-1},1\right\} \left(a^{\gamma}+b^{\gamma}\right)\quad\text{for all}\;a,b\ge0,\gamma>0,
\]
we continue to estimate (\ref{eq:2.9}) by
\begin{align}
\frac{1}{2^{\gamma}}\mathbf{A}^{\gamma} & \le C\left(\left\Vert \mathcal{A}_{y}^{\frac{k}{\gamma}}\mathbf{u}\left(a,\cdot\right)\right\Vert _{L^{2}\left(\Omega_{y}\right)}^{2}+\left\Vert \mathcal{A}_{y}^{\frac{k}{\gamma}-1}\mathbf{u}_{x}\left(a,\cdot\right)\right\Vert _{L^{2}\left(\Omega_{y}\right)}^{2}\right)^{\gamma}\nonumber \\
 & \le C\max\left\{ 2^{\gamma-1},1\right\} \left(\left\Vert \mathcal{A}_{y}^{\frac{k}{\gamma}}\mathbf{u}\left(a,\cdot\right)\right\Vert _{L^{2}\left(\Omega_{y}\right)}^{2\gamma}+\left\Vert \mathcal{A}_{y}^{\frac{k}{\gamma}-1}\mathbf{u}_{x}\left(a,\cdot\right)\right\Vert _{L^{2}\left(\Omega_{y}\right)}^{2\gamma}\right).\label{eq:2.10}
\end{align}

In case $k\equiv0\left(\text{mod}\;\gamma\right)$, the regularity
bound in (\ref{eq:2.10}) reduces to
\[
\mathbf{A}^{\gamma}\le2^{2\gamma}C\left(\left\Vert \mathbf{u}\left(a,\cdot\right)\right\Vert _{H^{\frac{k}{\gamma}}\left(\Omega_{y}\right)}^{2\gamma}+\left\Vert \mathbf{u}\left(a,\cdot\right)\right\Vert _{H^{\frac{k}{\gamma}-1}\left(\Omega_{y}\right)}^{2\gamma}\right).
\]

Similar to the case $k\le\gamma$, we deduce from (\ref{eq:2.9})
that
\[
\mathbf{A}^{\gamma}\le2^{2\gamma}C\left(\left\Vert \mathcal{A}_{y}^{\frac{k}{\gamma}}\mathbf{u}\left(a,\cdot\right)\right\Vert _{L^{2}\left(\Omega_{y}\right)}^{2\gamma}+\lambda_{1}^{\frac{k}{\gamma}-1}\left\Vert \mathbf{u}_{x}\left(a,\cdot\right)\right\Vert _{L^{2}\left(\Omega_{y}\right)}^{2\gamma}\right).
\]

Notice that when $k=0$ investigated in \cite{TTK15} , (\ref{eq:2.9})
gives us directly the natural criterion $\mathbf{u}\in C^{1}\left(\overline{\Omega_{x}};L^{2}\left(\Omega_{y}\right)\right)$.
\end{rem}

\section{Conclusions}

In general, solving the Cauchy problems of elliptic equations is doable
by regularization methods. In the context of kernel-based regularization,
this has been done in \cite{TTK15} and \cite{KTDT17}, working with
hardly checked criteria for convergence. In this Note, we have alleviated
such conditions, as informed in (\ref{eq:1.4})-(\ref{eq:1.5}) and
(\ref{eq:cond}), by the accepted regularity bound in Theorem \ref{thm:Consider-the-general}.
The ``accepted'' means that instead of testing the Gevrey-type criteria
on the exact solution, it now reduces to working with the forcing
function $f$. It therefore yields qualitatively better information
than previously developed assumptions. Interestingly, the Gevrey-type
criteria on $f$ can be ignored in computational environments by the
truncated Fourier series with the cut-off constant $N$ dependent
of the noise $\varepsilon$ assumed in (\ref{eq:epsilon}). Moreover,
the choice of $N$ can follow the work \cite{TTKT15}. Therefore,
one only needs $\mathbf{u}\in C^{1}\left(\overline{\Omega_{x}};H^{k}\left(\Omega_{y}\right)\right)\cap L^{1}\left(\Omega_{x};L^{2}\left(\Omega_{y}\right)\right)$
for any $k\in\mathbb{N}^{*}$ to solve the problem under consideration. 

As analyzed numerically in \cite{KTDT17}, this type of problems is
extremely sensitive to the noise level and the convergence is greatly
influenced by the boundedness of involved coefficients. It is worth
mentioning that the upper bound of the new criterion in Theorem \ref{thm:Consider-the-general}
still varies when doing with the truncated Fourier series on $f\left(\mathbf{u}\right)$.
Consequently, it may impact ugly on the theoretically desired convergence
of the proposed approximation. This unsolved issue will thus be our
next aim of study in the near future.
\begin{acknowledgement*}
This work is dedicated to the memory of V.A.K's father. The authors
desire to thank the handling editor and anonymous referee for their
helpful comments on this research.
\end{acknowledgement*}
\bibliographystyle{plain}
\bibliography{mybib}

\end{document}